\documentclass[12pt]{amsart}
\usepackage{graphicx, amssymb}
\numberwithin{equation}{section}

\newtheorem{Thm}{Theorem}[section]
\newtheorem{Rmk}{Remark}[section]
\newtheorem{Prop}{Proposition}[section]

\newtheorem{Lem}{Lemma}[section]

\def \N{{\mathbf N}}
\def \Q{{\mathbf Q}}
\def \R{{\mathbf R}}
\def \S{{\mathbf S}}
\def \T{{\mathbf T}}
\def \Z{{\mathbf Z}}

\def \w*L{w^*\hbox{-}L}

\def \a {{\alpha}}

\def \eps{{\epsilon}}

\def \th {\theta}
\def \O {\Omega}
\def \g {{\gamma}}

\def \l {{\lambda}}

\def \d{{\partial}}

\def\indc {{\bf 1}}

\begin{document}
\title[Free Path Lengths for Periodic Lorentz Gas]
      {On the Distribution of Free Path Lengths\\
       for the Periodic Lorentz Gas III}

\author[E. Caglioti]{Emanuele Caglioti}
\address[E. C.]%
{Dipartimento di Matematica\\
Istituto Guido Castelnuovo\\
Universit\`a di Roma ``La Sapienza"\\
p.le Aldo Moro 2, I00185 Roma}
\email{caglioti@mat.uniroma1.it}

\author[F. Golse]{Fran\c cois Golse}
\address[F. G.]%
{Institut Universitaire de France \\
 \& D\'epartement de Math\'e- matiques et Applications \\
 Ecole Normale Sup\'erieure Paris \\
 45 rue d'Ulm, F75230 Paris cedex 05}
\email{golse@dma.ens.fr}

\begin{abstract}
For $r\in(0,1)$, let $Z_r=\{x\in\R^2\,|\,\hbox{dist}(x,\Z^2)>r/2\}$
and $\tau_r(x,v)=\inf\{t>0\,|\,x+tv\in\d Z_r\}$. Let $\Phi_r(t)$ be
the probability that $\tau_r(x,v)\ge t$ for $x$ and $v$ uniformly
distributed in $Z_r$ and $\S^1$ respectively. We prove in this paper
that
$$
\begin{aligned}
\limsup_{\eps\to 0^+}
\frac1{|\ln\eps|}\int_\eps^{1/4}\Phi_r\left(\frac{t}{r}\right)\frac{dr}r
&=\frac2{\pi^2t}+O\left(\frac1{t^2}\right)
\\
\liminf_{\eps\to 0^+}
\frac1{|\ln\eps|}\int_\eps^{1/4}\Phi_r\left(\frac{t}{r}\right)\frac{dr}r
&=\frac2{\pi^2t}+O\left(\frac1{t^2}\right)
\end{aligned}
$$
as $t\to+\infty$. This result improves upon the bounds on $\Phi_r$ in
Bourgain-Golse-Wennberg [Commun. Math. Phys. 190 (1998), 491--508]. We
also discuss the applications of this result in the context of kinetic
theory.
\end{abstract}

\maketitle

\section{Statement of the problem and main results}\label{INTRO}

\subsection{The periodic Lorentz gas}

Let $r\in (0,\tfrac12 )$ and define
\begin{equation}
\label{Def-Z_r}
Z_r=\{x\in\R^2\,|\,\hbox{dist}(x,\Z^2)>r/2\}\,.
\end{equation}
Consider a point particle moving at speed $1$ inside $Z_r$ and
being specularly reflected each time it meets the boundary of
$Z_r$. Such a dynamical system is referred to as ``a periodic,
two-dimensional Lorentz gas". (Indeed, Lorentz used the methods
of kinetic theory to describe the motion of electrons in a metal
as that of a collisionless gas of point particles bouncing on
the crystalline structure of atoms in the metal \cite{Lo}).

The ``free path length" (or ``(forward) exit time") starting
from $x\in Z_r$ in the direction $v\in\S^1$'' is defined as
\begin{equation}
\label{Def-tau_r}
\tau_r(x,v)=\inf\{t>0\,|\,x+tv\in\d Z_r\}\,,
    \quad (x,v)\in Z_r\times\S^1\,.
\end{equation}
For each $v=(v_1,v_2)\in\S^1$ such that $v_1v_2\not=0$ and the
ratio $v_1/v_2\in\R\setminus\Q$, one has $\tau_r(x,v)<+\infty$,
since any orbit of a linear flow with irrational slope on the
2-torus is dense (see for instance in \cite{Ar}, --- section 51,
corollary 1, p. 287 --- for this well-known fact).

Set $Y_r=Z_r/\Z^2$; since $\tau_r(x,v)=\tau_r(x+k,v)$ for each
$(x,v)\in Z_r\times\S^1$ and $k\in\Z^2$, the function $\tau_r$
can be seen as defined on $Y_r\times\S^1$.

\begin{figure}
    \centering
    \includegraphics[height=8.0cm]{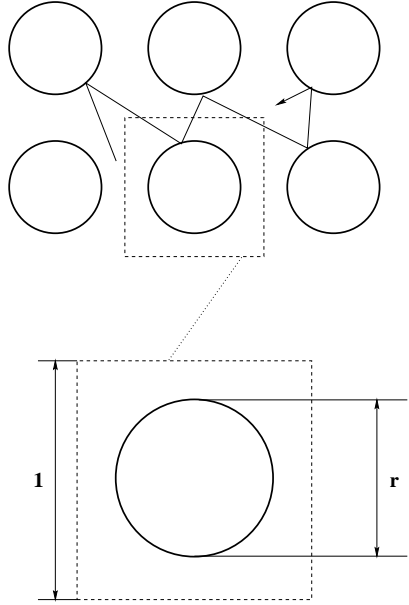}
    \caption{The Lorentz gas: $Z_r$ and the punctured torus $Y_r$}
\end{figure}

\subsection{Invariant measure for the Lorentz gas}

Let
$$
V_r=dxdv-\hbox{meas}(Y_r\times\S^1)
    \hbox{ and }\mu_r=\frac1{V_r}dxdv\,.
$$
Thus $\mu_r$ is a Borelian probability measure on $Y_r\times\S^1$.

On the other hand, the evolution of the Lorentz gas is governed
by the broken Hamiltonian flow
\begin{equation}
\label{Def-Flow}
\begin{aligned}
{}&(\dot{x}(t),\dot{v}(t))=(v(t),0)\,,
    &&\quad\hbox{ whenever }x(t)\in Y_r\,,
\\
&v(t^+)=v(t^-)-2v(t^-)\cdot n_{x(t)}n_{x(t)}\,,
    &&\quad\hbox{ whenever }x(t)\in\d Y_r\,,
\end{aligned}
\end{equation}
where $n_x$ denotes the inward unit normal to $\d Y_r$ at point $x$.
As can be easily checked, the measure $\mu_r$ is invariant under the
flow (\ref{Def-Flow}).

Let $m\in L^\infty(\S^1)$ such that
\begin{equation}
\label{cond-m}
m\ge 0\hbox{ and }\iint_{Y_r\times\S^1}m(\th)d\mu_r(x,\th)=1\,.
\end{equation}

The object of interest in the present paper is the distribution of
$\tau_r$ under $md\mu_r$, ie. the function $\Phi^m_r:\,\R_+\to[0,1]$
defined by
\begin{equation}
\label{Def-Phi_r}
\Phi^m_r(t)=md\mu_r-\hbox{Prob }
    (\{(x,v)\in Y_r\times\S^1\,|\,\tau_r(x,v)>t\})\,.
\end{equation}
More precisely, we are interested in the asymptotic behavior of
$\Phi^m_r(t)$ as $r\to 0^+$ and for large values of $t$ --- large
compared to $1/r$.

As explained in \cite{Bu} pp. 221--222, knowing the distribution of
the free path lengths in the small $r$ limit has several important
implications. For instance, it leads to the correct asymptotic model
for the Boltzmann-Grad limit of the periodic Lorentz gas --- which
is {\it not} governed by the linear Boltzmann equation (equation (10)
of \cite{Lo}): see below. Another application bears on the asymptotic
behavior of the Kolmogorov-Sinai entropy --- or equivalently, of the
Lyapunov exponent --- of the periodic Lorentz gas in the small $r$
limit. We refer to \cite{Ch} for a survey of the most recent results
and open questions on this subject, and to \cite{CG} which addresses
one of these open problems by methods similar to those developed here.

\subsection{Main result}

The reference \cite{BGW} established for each $m\in L^\infty(\S^1)$
as in (\ref{cond-m}) the existence of two positive constants $C_m$
and $C'_m$ such that, for each $r\in (0,\tfrac12)$ and each $t>1/r$
\begin{equation}
\label{Bound-phi_r}
\frac{C_m}{rt}\le\Phi^m_r(t)\le\frac{C'_m}{rt}\,;
\end{equation}
(see Theorems B and C in \cite{BGW}). In (\ref{Bound-phi_r}), the upper
bound was proved by an argument based on Fourier series, while the lower
bound was obtained by a geometric construction exhausting all possible
channels, ie. infinitely long open strips included in $Z_r$. The validity
of the lower bound in (\ref{Bound-phi_r}) was extended to arbitrary space
dimensions in \cite{GW}.

Numerical simulations in \cite{GW} suggest the following questions: for
each $m\in L^\infty(\S^1)$ as in (\ref{cond-m}) does one have, for each
$t$ large enough (say, for each $t>2$)
\begin{equation}
\label{Q1}
\Phi^m_r\left(\frac{t}r\right)\to\Lambda^m(t)\hbox{ as }\eps\to 0\,?
\end{equation}
and, if so, does one have, for some constant $C>0$
\begin{equation}
\label{Q2}
\Lambda^m(t)\sim\frac{C}t\hbox{ as }t\to+\infty\,?
\end{equation}

We have not been able to fully answer (\ref{Q1}), but, were (\ref{Q1})
true, our main result in this paper (Theorem \ref{LIM-DIST} below) would
answer (\ref{Q2}) by giving an explicit value for $C$. It confirms the
numerical results obtained in \cite{GW} (see figures 4 and 5 there).
Throughout the paper, the notation for the convergence in the sense of
Cesaro is as follows:
$$
C\!\!-\!\!\!\limsup_{\eps\to 0^+}f(\eps)=l\hbox{ means that }
\limsup_{\eps\to 0^+}\frac1{|\ln\eps|}\int_\eps^{\eps^*}f(r)\frac{dr}r=l
$$
for some $\eps^*>0$. (A similar notation is used for the $\liminf$ and the
$\lim$ in the sense of Cesaro).

\begin{Thm}\label{LIM-DIST}
Let $m\in L^\infty(\S^1)$ satisfying (\ref{cond-m}). Let $t^*>\sqrt{2}$ and let
$$
\Lambda^m_+(t^*)=C\!\!-\!\!\!\limsup_{r\to 0^+}\Phi^m_r\left(\frac{t^*}{r}\right)\,,
\hbox{ and }
\Lambda^m_-(t^*)=C\!\!-\!\!\!\liminf_{r\to 0^+}\Phi^m_r\left(\frac{t^*}{r}\right)\,.
$$
Then
\begin{equation}
\label{lim-dist}
\Lambda^m_+(t^*)\sim\frac2{\pi^2t^*}\,,\hbox{ and }
    \Lambda^m_-(t^*)\sim\frac2{\pi^2t^*}\hbox{ as }t^*\to+\infty\,.
\end{equation}
More precisely,
\begin{equation}
\label{asympt-dist}
\left|\Lambda^m_\pm(t^*)-\frac2{\pi^2t^*}\right|\le\frac{8\|m\|_{L^\infty}}{t^*-3}\,.
\end{equation}
\end{Thm}

A serious shortcoming of the result above is the need for averaging in $r$ before
letting $r\to 0^+$. It seems however that it cannot be avoided, at least by using
the techniques of the present paper.

This leads to a natural question, that of the choice of the measure $\frac{dr}r$
to define the Cesaro mean in Theorem \ref{LIM-DIST}. The reasons for this choice
are made clear by following the proof, but we take this opportunity of giving an
idea of this proof by providing some qualitative argument in favor of this choice.

The proof of Theorem \ref{LIM-DIST} is based on comparing the size $r$ of the
obstacle with the sequence of errors $d_n$ in the approximation by continued
fractions of $v_2/v_1$, ie. of the slope of the linear flow (say, in the case
where $0<v_2<v_1$). It is natural in this context to renormalize the problem
by applying to $\a=v_2/v_1$ the Gauss map $T:\,x\mapsto\frac1x-[\frac1x]$ ---
we refer to the appendix for more details on these notions. Lemma \ref{LEM-REN}
below and especially the formula $d_n(\a)=\a d_{n-1}(T(\a))$ show that the exit
time problem with slope $\a$ and obstacle of size $r$ is mapped to the analogous
problem with slope $T\a$ and obstacle of size $\a r$. Hence it is natural to
define the Cesaro average in Theorem \ref{LIM-DIST} with the measure $\frac{dr}r$
which is the scale invariant (Haar) measure of the multiplicative group $\R_+^*$.

\begin{Rmk}\label{CES-LIM}
In fact, one can prove that $\Lambda_+^m(t^*)=\Lambda_-^m(t^*)$ for each $t^*$
large enough (say, for $t^*>10$). (In other words $C\!\!-\!\!\!\lim_{r\to 0^+}
\Phi^m_r\left(\frac{t^*}{r}\right)$ exists). However, this result requires some
significant improvements of the method used in the present paper; they will be
described in \cite{CG}.
\end{Rmk}

\section{A partition of $\T^2$}

In 1989, R. Thom posed the following problem:
\begin{itemize}

\item{\it What is the longest orbit of a linear flow with irrational slope on a flat torus
with a disk removed?}

\end{itemize}
The answer to this question was found by Blank and Krikorian in \cite{BK} and is summarized
as follows. Without loss of generality, assume that the linear flow is $x\mapsto x+tv$ where
$v=(\cos\th,\sin\th)$ with $\th\in(0,\tfrac{\pi}{4})$. The removed disk of diameter $r$ is
replaced by a vertical slit $S_r(v)$ of length $r/\cos\th$ with the same center (see figure
2).

\begin{figure}
    \centering
    \includegraphics[height=8.0cm]{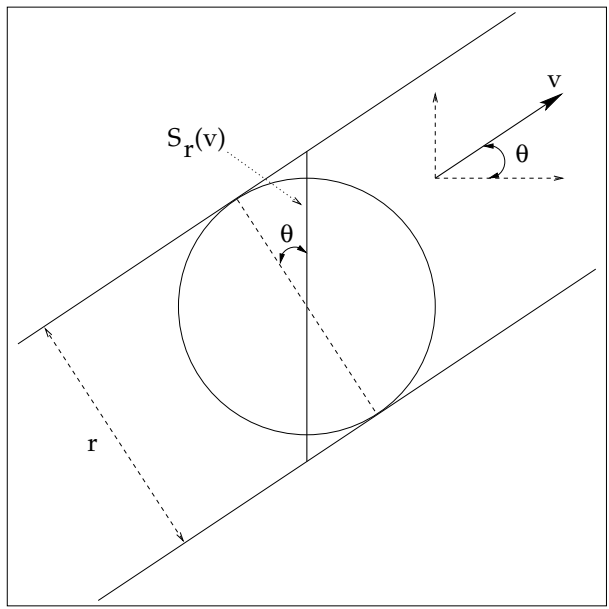}
    \caption{The punctured torus $Y_r$ and the slit $S_r(v)$}
\end{figure}

\begin{Prop}[Blank-Krikorian \cite{BK} p. 722]\label{BK-PROP}
Let $r\in(0,\tfrac12)$, $\th\in(0,\tfrac{\pi}{4})$ and $v=(\cos\th,\sin\th)$. Assume
that $\tan\th$ is irrational. There exist three positive numbers $l_A(r,v)$, $l_B(r,v)$
and $l_C(r,v)$ satisfying
$$
l_A(r,v)<l_B(r,v)\quad\hbox{and}\quad l_C(r,v)=l_A(r,v)+l_B(r,v)
$$
and such that, for any orbit $\g$ of the linear flow $x\mapsto x+tv$ in $\T^2\setminus S_r(v)$,
$\hbox{length}(\g)$ takes one of the three values $l_A(r,v)$, $l_B(r,v)$ or $l_C(r,v)$.

Conversely, given any $l\in\{l_A(r,v),l_B(r,v),l_C(r,v)\}$, there exists an orbit $\g$ of the
flow $x\mapsto x+tv$ of length $l$.
\end{Prop}

Let $v$ be fixed; orbits of length $l_A(r,v)$ (resp. of length $l_B(r,v)$, $l_C(r,v)$) are
referred to as orbits of type $A$ (resp. of type $B$, $C$)\footnote{On p. 722 of \cite{BK},
the sentence ``If the moving slit first meets the fixed slit at the bottom, the $A$ and
$B$-orbits are reversed" might be the source of a slight ambiguity in the definition of the
partition above. In the present paper, the orbits of type $A$ are the shortest, consistently
with the table on p. 726 of \cite{BK}. Hence the roles of orbits of type $A$ and $B$ cannot
be reversed in the present discussion.}. Proposition \ref{BK-PROP} defines a partition
$(Y_A(r,v),Y_B(r,v),Y_C(r,v))$ of $\T^2\setminus S_r(v)$, where
$$
Y_A(r,v)=\{x\in Y_r\setminus S_r(v)\,|\,x\hbox{ belongs to an orbit of type }A\}
$$
and $Y_B(r,v)$ and $Y_C(r,v)$ are similarly defined. Define further
$$
S_A(r,v)=\{y\in S_r(v)\,|\,\hbox{ the $v$-orbit starting from $y$ is of type A}\}
$$
with analogous definitions for $S_B(r,v)$ and $S_C(r,v)$. Clearly $Y_A(r,v)$ (resp. $Y_B(r,v)$,
$Y_C(r,v)$) is metrically equivalent to a strip (parallellogram) of length $l_A(r,v)$ and width
$|S_A(r,v)|$ (resp. of length $l_B(r,v)$, $l_C(r,v)$ and width $|S_B(r,v)|$, $|S_B(r,v)|$): see
figure 3 below.

\begin{figure}
    \centering
    \includegraphics[height=8.0cm]{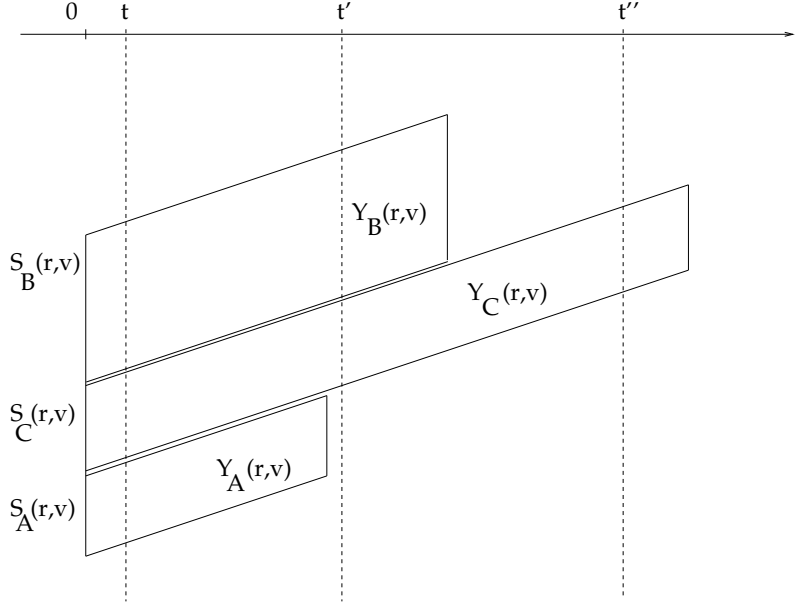}
    \caption{The partition $(Y_A(r,v),Y_B(r,v),Y_C(r,v))$ of $\T^2\setminus S_r(v)$}
\end{figure}

Define now $\l_r$ to be the exit time in the torus with the slit (instead of the disk)
removed,
ie.
\begin{equation}
\label{Def-l_r}
\l_r(z,v)=\inf\{t>0\,|\,z+tv\in S_r(v)\}\hbox{ for each }z\in\T^2\setminus S_r(v)\,.
\end{equation}
For $v=(\cos\th,\sin\th)$ with $\th\in(0,\tfrac{\pi}{4})$ as above, define
\begin{equation}
\label{DefPsi}
\psi_r(t,v)=dx-\hbox{Prob}(\{z\in\T^2\setminus S_r(v)\,|\,\l_r(z,v)\ge t\})\,.
\end{equation}
With the partition of $\T^2\setminus S_r(v)$ in $(Y_A(r,v),Y_B(r,v),Y_C(r,v))$ , which
is metrically equivalent to the disjoint union of three strips as represented on figure 3,
computing $\psi_r(t,v)$ in terms of the quantities $|S_A(r,v)|$, $l_A(r,v)$, $|S_B(r,v)|$,
$l_B(r,v)$, $|S_C(r,v)|$ and $l_C(r,v)$ becomes an easy task. One finds that

\begin{itemize}

\item if $0\le t\le l_A(r,v)$, then
\begin{equation}
\label{Psi1}
\psi_r(t,v)=1-tr\,;
\end{equation}

\item if $l_A(r,v)\le t'\le l_B(r,v)$, then
\begin{equation}
\label{Psi2}
\begin{aligned}
\psi_r(t',v)=1-&l_A(r,v)r
\\
-&(t'-l_A(r,v))(|S_B(r,v)|+|S_C(r,v)|)\cos\th\,;
\end{aligned}
\end{equation}

\item if $l_B(r,v)\le t''\le l_C(r,v)$, then
\begin{equation}
\label{Psi3}
\psi_r(t'',v)=(l_C(r,v)-t'')|S_C(r,v)|\cos\th\,;
\end{equation}

\item if $t\ge l_C(r,v)$, then
\begin{equation}
\label{Psi4}
\psi_r(t,v)=0\,.
\end{equation}

\end{itemize}

The graph of $t\mapsto\psi_r(t,v)$ is presented on figure 4 below.

\begin{figure}
    \centering
    \includegraphics[height=8.0cm]{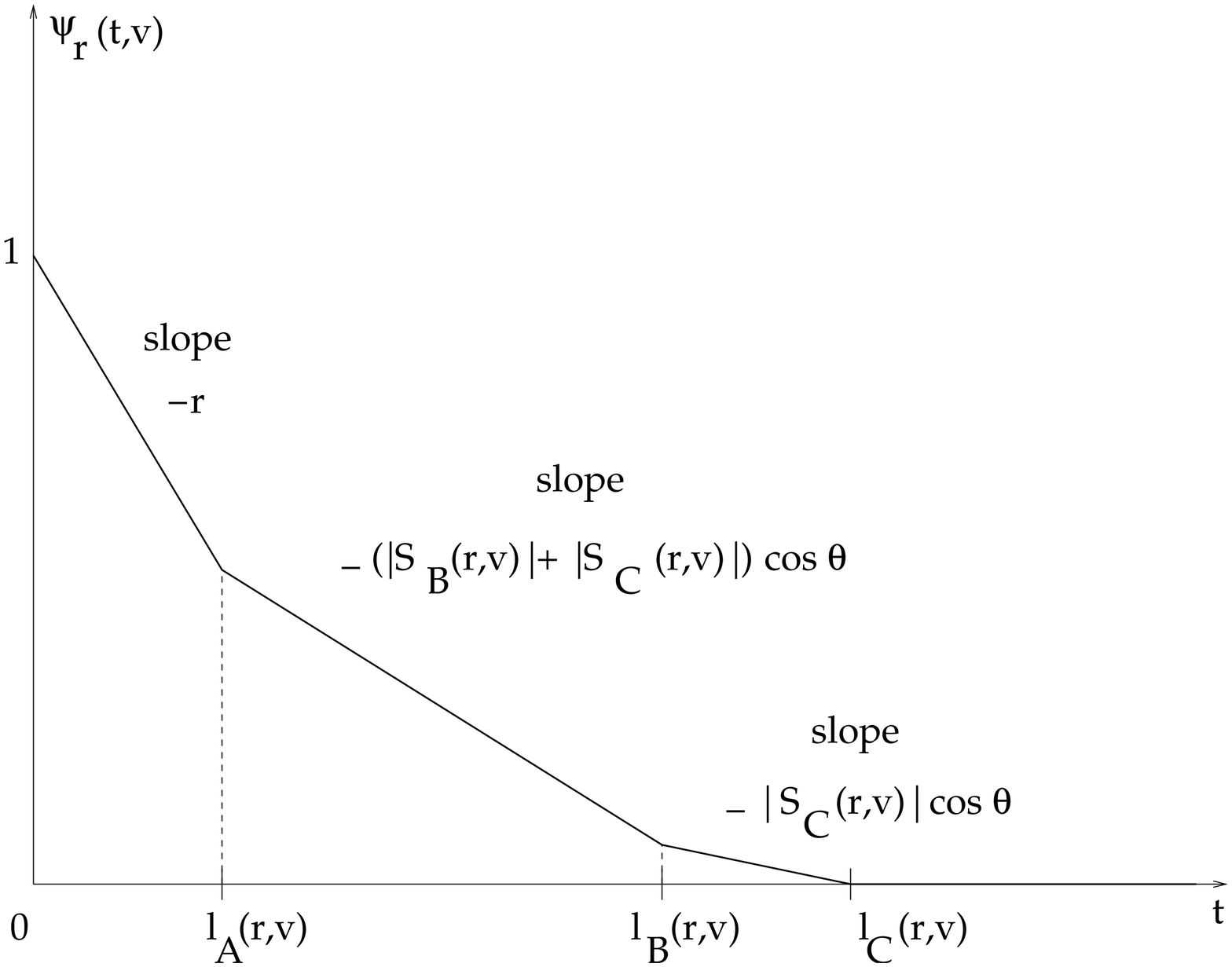}
    \caption{Graph of $t\mapsto\psi_r(t,v)$}
\end{figure}

So far, the distribution $\psi_r(t,v)$ has been computed in terms of the quantities
$|S_A(r,v)|$, $l_A(r,v)$, $|S_B(r,v)|$, $l_B(r,v)$, $|S_C(r,v)|$ and $l_C(r,v)$ that
characterize the partition $(Y_A(r,v),Y_B(r,v),Y_C(r,v))$ of $\T^2\setminus S_r(v)$.
These quantities can be expressed (see \cite{BK}) in terms of the continued fraction
expansion of $\tan\th$, in the following manner.

In the following discussion, we freely use the notations recalled in the Appendix
below. Let $\a=\tan\th\in(0,1)$; the $\th$'s for which $\a\in\Q$ form a set of
measure $0$ and are discarded in the argument below. With the sequence of errors
$d_n$ in the continued fraction expansion of $\a$, we consider the following
partition of the interval $(0,1)$:
\begin{equation}
\label{01-Part-a}
(0,1)=\bigcup_{n\ge 1}I_n\,,\quad\hbox{ with }I_n=[d_n,d_{n-1})\,.
\end{equation}
Each interval $I_n$ is further partitioned into
\begin{equation}
\label{01-Part-b}
I_n=\bigcup_{1\le k\le a_n}I_{n,k}\,,\quad\hbox{ with }
I_{n,k}=[\sup(d_n,d_{n-1}-kd_n),d_{n-1}-(k-1)d_n)\,,
\end{equation}
so that eventually we arrive at the following nested partition of $(0,1)$:
\begin{equation}
\label{01-Part}
(0,1)=\bigcup_{n\ge 1}\bigcup_{1\le k\le a_n}I_{n,k}\,.
\end{equation}

\begin{Prop}[Blank-Krikorian \cite{BK} p. 726]\label{BK-PROP2}
Let $\th\in(0,\tfrac{\pi}{4})$ and $r\in(0,\tfrac12)$. Let $v=(\cos\th,\sin\th)$ and
set $\a=\tan\th$ (it is assumed that $\a\notin\Q$). For $R=\frac{r}{\cos\th}$, the
integers $n\ge 1$ and $k$ such that $1\le k\le a_n$ ($a_n$ being the $n$-th term in
the continued fraction expansion of $\a$) are defined by the fact that $R\in I_{n,k}$
(since the intervals $I_{n,k}$ form a partition of $(0,1)$: see (\ref{01-Part})).
Then

\begin{itemize}

\item $l_A(r,v)=q_n$ and $|S_A(r,v)|=R-d_n$,
\item $l_B(r,v)=q_{n-1}+kq_n$ and $|S_B(r,v)|=R-(d_{n-1}-kd_n)$,
\item $l_C(r,v)=q_{n-1}+(k+1)q_n$ and $|S_C(r,v)|=d_{n-1}-(k-1)d_n-R$.

\end{itemize}
\end{Prop}

Using these values, we arrive at the following expression for $\psi_r(t,v)$, whenever
$R=\frac{r}{\cos\th}\in I_{n,k}$ (see figure 5 below):

\begin{itemize}

\item if $0\le t\cos\th\le q_n$, then
\begin{equation}
\label{Psi1'}
\psi_r(t,v)=1-tr\,,
\end{equation}
\item if $q_n\le t\cos\th\le q_{n-1}+kq_n$, then
\begin{equation}
\label{Psi2'}
\psi_r(t,v)=1-Rq_n-d_n(t\cos\th-q_n)\,,
\end{equation}
\item if $q_{n-1}+kq_n\le t\cos\th\le q_{n-1}+(k+1)q_n$, then
\begin{equation}
\label{Psi3'}
\begin{aligned}
\psi_r(t,v)=1-&Rq_n-d_n[q_{n-1}+(k-1)q_n]
\\
-&[R-(d_{n-1}-(k-1)d_n)](t\cos\th-q_{n-1}-kq_n)\,,
\end{aligned}
\end{equation}
\item if $t\cos\th\ge q_{n-1}+(k+1)q_n$, then
\begin{equation}
\label{Psi4'}
\psi_r(t,v)=0\,.
\end{equation}

\end{itemize}

\begin{figure}
    \centering
    \includegraphics[height=8.0cm]{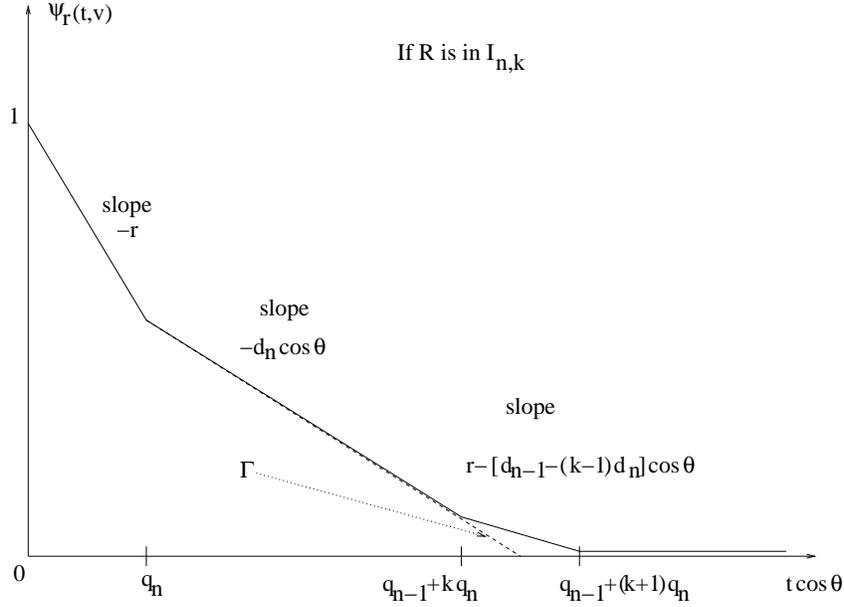}
    \caption{Graph of $t\mapsto\psi_r(t,v)$ for $R\in I_{n,k}$}
\end{figure}

\smallskip
Finally, the discussion above leads to the statistics of the exit time $\l_r(z,v)$
defined in (\ref{Def-l_r}) corresponding to the torus $\T^2$ with the slit $S_r(v)$
removed. In the case of the torus with the disk removed, the corresponding exit time
$\tau_r(x,v)$ defined in (\ref{Def-tau_r}) is related to $\l_r(z,v)$ by the obvious
inequalities
\begin{equation}
\label{tau-l-Ineq}
\l_r(x,v)-\tfrac{r}2\le\tau_r(x,v)\le\l_r(x,v)+\tfrac{r}2
\hbox{ for each }x\in Y_r\setminus S_r(v)\,.
\end{equation}
Define
\begin{equation}
\label{Def-phi_r}
\phi_r(t,v)=dx-\hbox{Prob}(\{x\in Y_r\setminus S_r(v)\,|\,\tau_r(x,v)\ge t\})\,;
\end{equation}
because of (\ref{tau-l-Ineq}), one has
\begin{equation}
\label{phi-psi-Ineq}
\psi_r(t+\tfrac{r}2,v)\le\phi_r(t,v)\le\psi_r(t-\tfrac{r}2,v)\,,\quad t\ge\tfrac{r}2\,.
\end{equation}
The remaining part of the paper uses the evaluation of $\phi_r$ based on this inequality
together with the formulas (\ref{Psi1'}), (\ref{Psi2'}), (\ref{Psi3'}) and (\ref{Psi4'})
for $\psi_r$.

\section{An ergodic theorem}\label{ERGO}

Given $\a\in(0,1)$ and $\eps\in(0,1)$, we define
\begin{equation}
\label{DefN}
N(\a,\eps)=\inf\{n\in\N\,|\,d_{n+1}(\a)<\eps\}\,.
\end{equation}
In terms of the partition (\ref{01-Part-a}) of $(0,1)$, $N(\a,\eps)$ can be equivalently
defined by the condition
\begin{equation}
\label{DefN'}
\eps\in I_{N(\a,\eps)+1}\,.
\end{equation}

\smallskip
We start by recalling the following more or less classical lemma.

\begin{Lem}\label{LEM-N}
For ae. $\a\in(0,1)$, one has
$$
N(\a,\eps)\sim -\tfrac{12\ln 2}{\pi^2}\ln\eps\quad\hbox{ as }\eps\to 0^+\,.
$$
\end{Lem}

\begin{proof}
The definition of $N(\a,\eps)$ and the third formula in Lemma \ref{LEM-REN} imply that
\begin{equation}
\label{IneqN}
\sum_{j=0}^{N(\a,\eps)-1}-\ln T^j\a\le -\ln\eps<\sum_{j=0}^{N(\a,\eps)}-\ln T^j\a\,,
\end{equation}
where $T$ is the Gauss map (\ref{DefGauss}). First we prove that $N(\a,\eps)\to+\infty$
as $\eps\to 0^+$ ae. in $\a\in(0,1)$. Indeed, let $C>0$ and let
$$
E_C=\{\a\in(0,1)\,|\,N(\a,\eps)\le C\,,\,\,\hbox{ for all }\eps\in (0,\tfrac12 )\}\,.
$$
For all $\a\in E_C$ and all $\eps\in(0,\tfrac12)$
$$
-\ln\eps\le\sum_{j=0}^C-\ln T^j\a\,;
$$
therefore, for all $\eps\in(0,\tfrac12)$,
$$
-\ln\eps\cdot dg-\hbox{meas }(E_C)\le(C+1)\int_0^1(-\ln\a)dg(\a)<+\infty\,,
$$
since the measure $dg$ in (\ref{InvMeas}) is invariant under $T$. This implies that
\begin{equation}
\label{ENegl}
\hbox{for each }C>0\,,\quad dg-\hbox{meas }(E_C)=0\,.
\end{equation}
Because $N(\a,\eps)$ is a nonincreasing function of $\eps$,
\begin{equation}
\label{NtoInfty}
\hbox{for each }\a\in\left(\bigcup_{m\ge 1}E_m\right)^c\,, N(\a,\eps)\to+\infty
\end{equation}
as $\eps\to 0^+$ and
\begin{equation}
\label{Ncvpp}
dg-\hbox{meas }\left(\bigcup_{m\ge 1}E_m\right)=0\,.
\end{equation}
Secondly, by Birkhoff's ergodic theorem, there exists a $dg$-negligeable set $E'$ such
that
$$
\hbox{for each }\a\in(E')^c\,,
        \quad\frac1N \sum_0^{N-1}-\ln T^j\a\to\int_0^1(-\ln\a)dg(\a)
$$
as $N\to+\infty$, since the Gauss transformation $T$ of $(0,1)$ is ergodic with respect
to the invariant measure $dg(\a)$ --- see the Appendix below. Hence
$$
\begin{aligned}
\hbox{for each }\a\in\left(E'\cup\bigcup_{m\ge 1}E_m\right)^c&\,,
\\
\frac1{N(\a,\eps)}\sum_0^{N(\a,\eps)-1}-\ln T^j\a&\to\int_0^1(-\ln\a)dg(\a)\,.
\end{aligned}
$$
By (\ref{IneqN}), (\ref{NtoInfty}) and (\ref{Ncvpp}), one finally obtains that, as
$\eps\to 0^+$,
$$
\hbox{for each }\a\in\left(E'\cup\bigcup_{m\ge 1}E_m\right)^c\,,
    \quad\frac{-\ln\eps}{N(\a,\eps)}\to\int_0^1(-\ln\a)dg(\a)=\frac{\zeta(2)}{2\ln 2}
$$
(replacing the $\ln$ under the integral sign by its Taylor series at $\a=1$).
\end{proof}

\smallskip
The main result of this section is the following application of the Birkhoff ergodic
theorem. We shall need the notations below:
\begin{equation}
\label{DefDlt}
\begin{aligned}
\Delta_0(\a,x)&=-x-\ln d_{N(\a,e^{-x})+1}(\a)\,,\quad
\\
\Delta_1(\a,x)&=-x-\ln d_{N(\a,e^{-x})}(\a)\,.
\end{aligned}
\end{equation}

\begin{Prop}\label{PROP-ERG}
Let $f$ be a bounded nonnegative measurable function on $\R^2$. For each $x^*\in\R$
and ae. in $\a\in(0,1)$, one has
$$
\frac1{|\ln\eps|}\int_{x^*}^{|\ln\eps|} f(\Delta_0(\a,x),\Delta_1(\a,x))dx
    \to\tfrac{12}{\pi^2}\int_0^1\frac{F(\th)d\th}{1+\th}
$$
as $\eps\to 0$, where
$$
F(\th)=\int_{0}^{|\ln(\th)|}f(|\ln(\th)|-y,-y)dy\,.
$$
\end{Prop}

\begin{proof}
First we decompose the integral
\begin{equation}
\label{IntDec}
\begin{aligned}
\int_{x^*}^{|\ln\eps|}f(\Delta_0(\a,x),&\Delta_1(\a,x))dx
\\
&=\int_{x^*}^0f(\Delta_0(\a,x),\Delta_1(\a,x))dx
\\
&+\int_{|\ln d_{N(\a,\eps)}(\a)|}^{|\ln\eps|}f(\Delta_0(\a,x),\Delta_1(\a,x))dx
\\
    &+\sum_{l=0}^{N(\a,\eps)-1}
    \int_{|\ln d_l(\a)|}^{|\ln d_{l+1}(\a)|}f(\Delta_0(\a,x),\Delta_1(\a,x))dx
\end{aligned}
\end{equation}
Now observe that, whenever $x$ belongs to the domain of integration of
$$
\int_{|\ln d_l(\a)|}^{|\ln d_{l+1}(\a)|}f(\Delta_0(\a,x),\Delta_1(\a,x))dx\,,
$$
in other words, whenever
$$
|\ln d_l(\a)|<x\le|\ln d_{l+1}(\a)|\,,\quad\hbox{ then }N(\a,e^{-x})=l\,.
$$
This implies that
\begin{equation}
\label{ChgVar}
\begin{aligned}
\int_{|\ln d_l(\a)|}^{|\ln d_{l+1}(\a)|}&f(\Delta_0(\a,x),\Delta_1(\a,x))dx
\\
&=
\int_{-\ln d_l(\a)}^{-\ln d_{l+1}(\a)}f(-x-\ln d_{l+1}(\a),-x-\ln d_l(\a))dx
\\
&=
\int_0^{\ln(d_l(\a)/d_{l+1}(\a))}f(\ln(d_l(\a)/d_{l+1}(\a))-y,-y)dy
\\
&=
F(T^l\th)\,.
\end{aligned}
\end{equation}
Likewise, since $f$ is nonnegative
\begin{equation}
\label{ChgVar2}
\begin{aligned}
0\le\int_{|\ln d_{N(\a,\eps)}(\a)|}^{|\ln\eps|}f(\Delta_0(\a,x),&\Delta_1(\a,x))dx
\\
&\le
\int_{|\ln d_{N(\a,\eps)}(\a)|}^{|\ln d_{N(\a,\eps)+1}(\a)|}
    f(\Delta_0(\a,x),\Delta_1(\a,x))dx
\\
&\le
F(T^{N(\a,\eps)}\a)\,.
\end{aligned}
\end{equation}
Using (\ref{IntDec}) and (\ref{ChgVar}) leads to
\begin{equation}
\label{IntDec2}
\begin{aligned}
\frac1{|\ln\eps|}\int_{x^*}^{|\ln\eps|}f(\Delta_0(\a,x),&\Delta_1(\a,x))dx
\\
&=\frac{N(\a,\eps)}{|\ln\eps|}
    \left(\frac1{N(\a,\eps)}\sum_{l=0}^{N(\a,\eps)}F(T^l\th)\right)
\\
&=\frac1{|\ln\eps|}\int_{x^*}^0f(\Delta_0(\a,x),\Delta_1(\a,x))dx
\\
&+\frac1{|\ln\eps|}\int_{|\ln d_{N(\a,\eps)}(\a)|}^{|\ln\eps|}
    f(\Delta_0(\a,x),\Delta_1(\a,x))dx\,.
\end{aligned}
\end{equation}
Observe first that, by its definition
$$
0\le F(\th)\le\|f\|_{L^\infty}|\ln\th|\hbox{ which implies that }F\in L^1((0,1),dg)\,.
$$
Birkhoff's ergodic theorem applied to the Gauss transformation $T$ together with
Lemma \ref{LEM-N} shows that
$$
\frac1{N(\a,\eps)}\sum_{l=0}^{N(\a,\eps)}F(T^l\th)
    \to\frac{1}{\ln 2}\int_0^1F(\th)dg(\th)
$$
for ae. $\a\in(0,1)$ as $\eps\to 0$, so that the first term in the right hand side
of (\ref{IntDec2}) converges to
$$
\tfrac{12}{\pi^2}\int_0^1F(\th)dg(\th)
$$
as $\eps\to 0^+$ ae. in $\a$. The second term in the right hand side of (\ref{IntDec2})
obviously vanishes ae. in $\a\in(0,1)$ as $\eps\to 0$. As for the third term, because
of (\ref{ChgVar2}), one has
$$
\begin{aligned}
0\le \frac1{|\ln\eps|}\int_{|\ln d_{N(\a,\eps)}(\a)|}^{|\ln\eps|}
    f(\Delta_0(\a,x),&\Delta_1(\a,x))dx
\\
&\le
\frac{N(\a,\eps)}{|\ln\eps|}\,\frac1{N(\a,\eps)}F(T^{N(\a,\eps)}\a)
\end{aligned}
$$
and the right hand side of the above inequality converges to 0 ae. in $\a\in(0,1)$ as
$\eps\to 0^+$ by Birkoff's ergodic theorem and Lemma \ref{LEM-N}.
\end{proof}

\section{Proof of theorem \ref{LIM-DIST}}\label{PROOF}

\subsection{Step 1: pointwise estimates}

The discussion in section \ref{ERGO} made it clear that the natural objects for applying
Birkhoff's ergodic theorem are functions involving a finite number (two in this case) of
the $d_n$'s with $n\to+\infty$. In this first step, we shall reduce the function $\psi_r$
to an expression of this form modulo terms that are small in some appropriate sense in the
asymptotic regime that we consider --- ie. as $t^*\to+\infty$.

Throughout this subsection, we set $r\in(0,\tfrac12)$ and $\th\in(0,\tfrac{\pi}4)$ such
that $\a=\tan\th\notin\Q$; let then $v=(\cos\th,\sin\th)$ and $R=\frac{r}{\cos\th}$. We
also consider $n\ge 1$ and $k\in\{1,\ldots,a_n\}$ (where $\a=[a_1,a_2\ldots]$) such that
$R\in I_{n,k}$ --- this defines $n$ and $k$ in a unique way, since the $I_{n,k}$'s form
a partition of $(0,1)$.

\begin{Lem}\label{INEQ-qn}
Under these conditions on $\a$, $n$ $k$ and $R$

\begin{itemize}

\item the integer $k$ is given by the formula
\begin{equation}
\label{Def-k}
k=\inf\{l\in\N^*\,|\,d_{n-1}-ld_n\le R\}=\left[\frac{d_{n-1}-R}{d_n}\right]+1\,;
\end{equation}

\item the denominator $q_n$ of the $n$-th convergent of $\a$ satisfies the estimate
\begin{equation}
\label{Estim-qn}
\frac1{R+(k+1)d_n}<q_n<\frac1{R+(k-1)d_n}\,;
\end{equation}

\item finally, one has
\begin{equation}
\label{Estim-qndn}
q_nd_n<\frac1k\,,\quad 0<1-q_nd_{n-1}<\frac2{k+1}\,.
\end{equation}
\end{itemize}
\end{Lem}

\begin{proof}
Whenever $d_n\le R<d_{n-1}$ --- ie. whenever $R\in I_n$ --- the condition $R\in I_{n,k}$
amounts to defining $k$ by the formula
$$
k=\inf\{l\in\N^*\,|\,d_{n-1}-ld_n\le R\}=\left[\frac{d_{n-1}-R}{d_n}\right]+1\,.
$$

Next we estimate $q_n$. Since $R\in I_{n,k}$, one has in particular the inequalities
\begin{equation}
\label{IneqR}
d_{n-1}-kd_n\le R<d_{n-1}-d_n(k-1)\,,
\end{equation}
which imply
$$
\frac1{R+d_n(k-1)}>\frac1{d_{n-1}}>q_n\,,
$$
by the second inequality in (\ref{IneqErr}). This is exactly the upper bound in
(\ref{Estim-qn}). By (\ref{qdRel}),
$$
\begin{aligned}
1-Rq_n-d_n(k-1)q_n&=d_nq_{n-1}+d_{n-1}q_n-Rq_n-d_n(k-1)q_n
\\
&<
q_n(d_n+d_{n-1}-R-(k-1)d_n)
\\
&<
q_n(d_n+d_{n-1}-d_{n-1}+kd_n-(k-1)d_n)
\\
&= 2d_nq_n
\end{aligned}
$$
which gives the lower bound in (\ref{Estim-qn}).

The upper bound in (\ref{Estim-qn}) and the fact that $R\ge d_n$ (since $R\in I_{n,k}$)
imply the first inequality in (\ref{Estim-qndn}). As for the second inequality there,
observe that $q_nd_{n-1}\le 1$ by (\ref{IneqErr}) which establishes the lower bound,
while the upper bound follows from (\ref{Estim-qn}) in the following manner:
$$
\begin{aligned}
\frac1{d_{n-1}}-q_n&\le\frac{R+d_n(k+1)-d_{n-1}}{d_{n-1}(R+d_n(k+1))}
\\
&<\frac{2d_n}{d_{n-1}(R+d_n(k+1))}<\frac{2}{d_{n-1}(k+1)}
\end{aligned}
$$
(where the penultimate inequality above follows from (\ref{IneqR})).
\end{proof}

\begin{Lem}\label{FORM-psir}
Let $v$, $\a$, $n$, $k$ and $R$ be chosen as above, and let $t^*>2$. Define
\begin{equation}
\label{def-chir}
\chi_r\left(\frac{t^*}{r},v\right)
=\left(1-\frac{R}{d_{n-1}}-t^*\frac{d_n}{R}\right)_+\,;
\end{equation}
then for $\frac{t^*}{R}\ge q_n$ we have
\begin{equation}
\label{psi-chir}
\left|\psi_r\left(\frac{t^*}{r},v\right)-\chi_r\left(\frac{t^*}{r},v\right)\right|
\le\frac4k \indc_{k\ge t^*-2}\,.
\end{equation}
\end{Lem}

\begin{proof}
Let $t^*>2$. In view of the second inequality in (\ref{IneqErr})
\begin{equation}
\label{Ineq1t*}
R<d_{n-1}<\frac1{q_n}\,,\quad\hbox{ so that }q_n<\frac1R<\frac{t^*}{R}\,.
\end{equation}
Next we compare the part of the graph of $t\mapsto\psi_r(t,v)$ that corresponds to
$t\cos\th\ge q_n$ with the straight line $\Gamma$ as defined on figure 5. They only
differ when $q_{n-1}+kq_n<t\cos\th<q_{n-1}+(k+1)q_n$, and since $\psi_r$ is a non
increasing function, they differ by at most
$$
\psi_r(\frac{q_{n-1}+kq_n}{\cos\th},v)
=
\left((d_{n-1}-(k-1)d_n)-R \right)q_n\,.
$$
By (\ref{IneqR}),
$$
0<[d_{n-1}-(k-1)d_n]-R<d_n\,.
$$
Hence
\begin{equation}
\label{Ineq2t*}
\begin{aligned}
0&\le\psi_r\left(\frac{t^*}{r},v\right)
    -\left(1-Rq_n-\left(\frac{t^*}{R}-q_n\right)d_n\right)_+
\\
&\le
\psi_r\left(\frac{q_{n-1}+kq_n}{\cos\th},v\right)
    \indc_{q_{n-1}+kq_n<\frac{t^*}{R}<q_{n-1}+(k+1)q_n}
\\
&\le
q_nd_n\indc_{t^*<(k+2)q_nR}\le\frac1k\indc_{k\ge t^*-2}\,,
\end{aligned}
\end{equation}
where the last inequality follows from (\ref{Estim-qndn}) and (\ref{Ineq1t*}).

Next we estimate the difference
$$
\left(1-Rq_n-\left(\frac{t^*}{R}-q_n\right)d_n\right)_+
-
\left(1-\frac{R}{d_{n-1}}-t^*\frac{d_n}{R}\right)_+
$$
Because of the second inequality in (\ref{IneqErr}), this difference is nonnegative.
Because the map $x\mapsto x_+$ is a contraction, this difference is less than
$$
\frac{R}{d_{n-1}}-Rq_n+q_nd_n\,;
$$
on the other hand both terms in the difference above vanish whenever $\frac{t^*}{R}
>q_{n-1}+(k+1)q_n$. Hence
$$
\begin{aligned}
\left|\left(1-Rq_n-\left(\frac{t^*}{R}-q_n\right)d_n\right)_+
-
\left(1-\frac{R}{d_{n-1}}-t^*\frac{d_n}{R}\right)_+\right|&
\\
\le
\left(\frac{R}{d_{n-1}}-Rq_n+q_nd_n\right)\indc_{\frac{t^*}{R}\le q_{n-1}+(k+1)q_n}&
\\
\le
\left(\frac2{k+1}+\frac1k\right)\indc_{\frac{t^*}{R}\le q_{n-1}+(k+1)q_n}
\le\frac3k \indc_{k\ge t^*-2}&
\end{aligned}
$$
where the penultimate inequality rests on (\ref{Estim-qndn}) and the fact that $R<d_{n-1}$.
\end{proof}

\subsection{Step 2: applying the ergodic theorem}

In this subsection again, we set $r\in(0,\tfrac12)$ and pick $\th\in(0,\tfrac{\pi}4)$ such
that $\a=\tan\th\notin\Q$; again we set $v=(\cos\th,\sin\th)$ and $R=\frac{r}{\cos\th}$. As
in the previous subsection, we define $n\ge 1$ by the condition $R\in I_n$. In other words,
we set $n=N(\a,R)$. With these assumptions, consider the expression $\chi_r(\frac{t^*}R,v)$
given by (\ref{def-chir}), ie.
$$
\begin{aligned}
\chi_r\left(\frac{t^*}R,v\right)&=
    \left(1-\frac{R}{d_{N(\a,R)-1}}-t^*\frac{d_{N(\a,R)}}{R}\right)_+
\\
&=
\left(1-e^{\Delta_1(\a,e^{-x})}-t^*e^{-\Delta_0(\a,e^{-x})}\right)_+
\end{aligned}
$$
with $x=-\ln R$ while $\Delta_0(\a,e^{-x})$ and $\Delta_1(\a,e^{-x})$ are defined as in
(\ref{DefDlt}).

\begin{Prop}\label{PROP-ERG2}
Let $R^*\in(0,1)$. Let $t^*>1$; then, for ae. $\th\in(0,\tfrac{\pi}4)$ such that $\a=
\tan\th\notin\Q$
$$
\begin{aligned}
\frac1{|\ln\eps|}\int_{\eps}^{R^*}&\chi_{R\cos\th}\left(\frac{t^*}R,v\right)\frac{dR}{R}
\to
\\
&\tfrac{12}{\pi^2}\int_0^1\left(\ln\frac{1+\sqrt{1-z}}{1-\sqrt{1-z}}
    -\sqrt{1-z}\right)\frac{dz}{4t^*+z}
\\
&+
\tfrac{6}{\pi^2}\int_0^1\left(\frac{z}{1+\sqrt{1-z}}
    -\frac{z}{1-\sqrt{1-z}}\right)\frac{dz}{4t^*+z}\,,
\end{aligned}
$$
as $\eps\to 0^+$, where $v=(\cos\th,\sin\th)$.
\end{Prop}

\begin{proof}
The proof is based upon applying Proposition \ref{PROP-ERG} to the function $f$ defined
by $f(z_1,z_2)=(1-e^{z_2}-t^*e^{-z_1})_+$, since, for each $R^*\in(0,1)$
$$
\begin{aligned}
\frac1{|\ln\eps|}\int_{\eps}^{R^*}\chi_{R\cos\th}&\left(\frac{t^*}R,v\right)\frac{dR}{R}
\\
&=
\frac1{|\ln\eps|}\int_{x^*}^{|\ln\eps|}f(\Delta_0(\a,e^{-x}),\Delta_1(\a,e^{-x}))dx\,.
\end{aligned}
$$
Starting from $f$, an elementary computation leads to $F$ defined for each $\xi\in(0,1)$
as in Proposition \ref{PROP-ERG} by
$$
\begin{aligned}
F(\xi)=\int_0^{|\ln\xi|}f(|\ln\xi|-y,-y)dy&=\int_0^{|\ln\xi|}(1-e^{-y}-t^*\xi e^y)_+dy
\\
&=
\int_1^{1/\xi}(\zeta-1-t^*\xi\zeta^2)_+\frac{d\zeta}{\zeta^2}\,.
\end{aligned}
$$
Assuming that $t^*>1$ and $\xi\in(0,1)$, elementary computations show that
$$
\begin{aligned}
\zeta-1-t^*\xi\zeta^2&<0
    \hbox{ for all }\zeta\in\R\hbox{ if }4t^*\xi>1\,,\hbox{ otherwise }
\\
\zeta-1-t^*\xi\zeta^2&\ge 0\hbox{ iff }
    \frac{1-\sqrt{1-4t^*\xi}}{2t^*\xi}\le\zeta\le\frac{1+\sqrt{1-4t^*\xi}}{2t^*\xi}
\end{aligned}
$$
and hence
$$
\begin{aligned}
F(\xi)&=\int_1^{1/\xi}(\zeta-1-t^*\xi\zeta^2)_+\frac{d\zeta}{\zeta^2}
\\
&=
\int_{\frac{1-\sqrt{1-4t^*\xi}}{2t^*\xi}}^{\frac{1+\sqrt{1-4t^*\xi}}{2t^*\xi}}
    \left(\frac1\zeta-\frac1{\zeta^2}-t^*\xi\right)d\zeta
\\
&=
\ln\frac{1+\sqrt{1-4t^*\xi}}{1-\sqrt{1-4t^*\xi}}-\sqrt{1-4t^*\xi}
\\
&+
\frac{2t^*\xi}{1+\sqrt{1-4t^*\xi}}-\frac{2t^*\xi}{1-\sqrt{1-4t^*\xi}}\,.
\end{aligned}
$$
Therefore
$$
\begin{aligned}
{}&\tfrac{12}{\pi^2}\int_0^1\frac{F(\xi)d\xi}{1+\xi}
\\
&=
\tfrac{12}{\pi^2}\int_0^{1/4t^*}\left(\ln\frac{1+\sqrt{1-4t^*\xi}}{1-\sqrt{1-4t^*\xi}}
    -\sqrt{1-4t^*\xi}\right)\frac{d\xi}{1+\xi}
\\
&+
\tfrac{12}{\pi^2}\int_0^{1/4t^*}\left(\frac{2t^*\xi}{1+\sqrt{1-4t^*\xi}}
    -\frac{2t^*\xi}{1-\sqrt{1-4t^*\xi}}\right)\frac{d\xi}{1+\xi}
\\
&=
\tfrac{12}{\pi^2}\int_0^1\left(\ln\frac{1+\sqrt{1-z}}{1-\sqrt{1-z}}
    -\sqrt{1-z}\right)\frac{dz}{4t^*+z}
\\
&+
\tfrac{6}{\pi^2}\int_0^1\left(\frac{z}{1+\sqrt{1-z}}
    -\frac{z}{1-\sqrt{1-z}}\right)\frac{dz}{4t^*+z}
\end{aligned}
$$
\end{proof}

\subsection{Step 3: $L^1$ estimate of the remainder}

The last ingredient in the proof of Theorem \ref{LIM-DIST} consists in estimating the
right hand side of (\ref{psi-chir}) in average (integrating over the angle $\th$).

Let $R\in(0,1)$ and $\a\in(0,1)\setminus\Q$; we define $n$ and $k$ by the condition
$R\in I_{n,k}$ --- below, this value of $k$ is denoted by $k(\a,R)$. Equivalently,
$$
n=N(\a,R)+1\,,\quad
    k(\a,R)=\left[\frac{d_{N(\a,R)}(\a)-R}{d_{N(\a,R)+1}(\a)}\right]+1\,.
$$

\begin{Lem}\label{LEM-remain}
Let $m\equiv m(\th)\in L^\infty((0,\tfrac{\pi}4))$. Then, for each $R\in(0,1)$ and
each $\l>1$, one has
\begin{equation}
\label{estim-remain}
\int_0^{\pi/4}\indc_{k(\tan\th,R)>\l}\,m(\th)d\th\le\frac{2\|m\|_{L^\infty}}{\l-1}\,.
\end{equation}
\end{Lem}

\begin{proof}
The definition of $k(\a,R)$ implies in particular that
$$
k(\a,R)\le\frac{d_{N(\a,R)}(\a)}{d_{N(\a,R)+1}(\a)}+1
    =\frac1{T^{N(\a,R)}\a}+1\,.
$$
Hence, for each $\l>1$, one has
$$
\begin{aligned}
dg-\hbox{meas }&(\{\a\in(0,1)\,|\,k(\a,R)\ge\l\})
\\
&\le
dg-\hbox{meas }\left(\left\{\a\in(0,1)\,|\,0<T^{N(\a,R)}\a\le\frac1{\l-1}\right\}\right)
\\
&\le
dg-\hbox{meas }\left(\left(0,\frac1{\l-1}\right]\right)
=
\frac1{\ln2}\ln\left(\frac\l{\l-1}\right)\to 0
\end{aligned}
$$
as $\l\to+\infty$. Changing variables from $\a$ to $\th=\arctan\a$ and using the classical
inequality $\ln(1+z)\le z$ leads to (\ref{estim-remain}).
\end{proof}

\subsection{Step 4: end of the proof}

We conclude the proof of Theorem \ref{LIM-DIST} by bringing together the various ingredients
described above. Let $m\in L^\infty([0,\tfrac{\pi}4])$ such that $m\ge 0$ and $\int_0^{\pi/4}
m(\th)d\th=1$.

By Lemma \ref{FORM-psir} --- and especially the inequality (\ref{psi-chir}) there ---  and
Lemma \ref{LEM-remain}, one has
$$
\begin{aligned}
\left|\frac1{\ln\eps}\int_\eps^{\eps^*}\int_0^{\pi/4}\right.&
    \psi_r\left(\frac{t^*}r,(\cos\th,\sin\th)\right)m(\th)d\th\frac{dr}{r}
\\
&\left.-
\frac1{\ln\eps}\int_\eps^{\eps^*}\int_0^{\pi/4}
    \chi_r\left(\frac{t^*}r,(\cos\th,\sin\th)\right)m(\th)d\th\frac{dr}{r}\right|
\\
&\le
\frac4{t^*-2}\,\frac1{|\ln\eps|}\int_\eps^{\eps^*}\int_0^{\pi/4}
    \indc_{k(\tan\th,r/\cos\th)>t^*-2}m(\th)d\th\frac{dr}{r}
\\
&\le\frac{8\|m\|_{L^\infty}}{(t^*-3)^2}
\end{aligned}
$$
By Proposition \ref{PROP-ERG2} and dominated convergence,
$$
\frac1{|\ln\eps|}\int_{\eps}^{\eps^*}\int_0^{\pi/4}
    \chi_r\left(\frac{t^*}r,(\cos\th,\sin\th)\right)m(\th)d\th\frac{dr}{r}
\to\Lambda(t^*)
$$
where
$$
\begin{aligned}
\Lambda(t^*)=&\tfrac{12}{\pi^2}\int_0^1\left(\ln\frac{1+\sqrt{1-z}}{1-\sqrt{1-z}}
    -\sqrt{1-z}\right)\frac{dz}{4t^*+z}
\\
+&
\tfrac{6}{\pi^2}\int_0^1\left(\frac{z}{1+\sqrt{1-z}}
    -\frac{z}{1-\sqrt{1-z}}\right)\frac{dz}{4t^*+z}\,,
\end{aligned}
$$
as $\eps\to 0^+$. Hence
$$
\limsup_{\eps\to 0^+}\frac1{|\ln\eps|}\int_{\eps}^{\eps^*}\int_0^{\pi/4}
    \psi_r\left(\frac{t^*}r,(\cos\th,\sin\th)\right)m(\th)d\th\frac{dr}{r}
\le
\Lambda(t^*)+\frac{8\|m\|_{L^\infty}}{(t^*-3)^2}\,,
$$
while
$$
\liminf_{\eps\to 0^+}\frac1{|\ln\eps|}\int_{\eps}^{\eps^*}\int_0^{\pi/4}
    \psi_r\left(\frac{t^*}r,(\cos\th,\sin\th)\right)m(\th)d\th\frac{dr}{r}
\ge
\Lambda(t^*)-\frac{8\|m\|_{L^\infty}}{(t^*-3)^2}\,.
$$
As $t^*\to+\infty$, one has
$$
\begin{aligned}
\Lambda(t^*)\sim&\frac{3}{\pi^2t^*}\int_0^1\left(\ln\frac{1+\sqrt{1-z}}{1-\sqrt{1-z}}
    -\sqrt{1-z}\right)dz
\\
+&
\frac{3}{2\pi^2t^*} \int_0^1\left(\frac{z}{1+\sqrt{1-z}}-\frac{z}{1-\sqrt{1-z}}\right)dz
=
\frac2{\pi^2t^*}\,.
\end{aligned}
$$
Since the function $t\mapsto\psi_r(t,v)$ is nonincreasing for all $v\in S^1$ and
$r\in(0,\tfrac12)$, one has, by using (\ref{phi-psi-Ineq}),
$$
\begin{aligned}
{}&\limsup_{\eps\to 0^+}\frac1{|\ln\eps|}\int_{\eps}^{\eps^*}\int_0^{\pi/4}
    \phi_r\left(\frac{t^*}r,(\cos\th,\sin\th)\right)m(\th)d\th\frac{dr}{r}
\\
=&\limsup_{\eps\to 0^+}\frac1{|\ln\eps|}\int_{\eps}^{2s}\int_0^{\pi/4}
    \phi_r\left(\frac{t^*}r,(\cos\th,\sin\th)\right)m(\th)d\th\frac{dr}{r}
\\
&\le
\limsup_{\eps\to 0^+}\frac1{|\ln\eps|}\int_{\eps}^{2s}\int_0^{\pi/4}
\psi_r\left(\frac{t^*-\tfrac{r}2}r,(\cos\th,\sin\th)\right)m(\th)d\th\frac{dr}{r}
\\
&\le
\limsup_{\eps\to 0^+}\frac1{|\ln\eps|}\int_{\eps}^{2s}\int_0^{\pi/4}
    \psi_r\left(\frac{t^*-s}r,(\cos\th,\sin\th)\right)m(\th)d\th\frac{dr}{r}
\\
&\le
\Lambda(t^*-s)+\frac{8\|m\|_{L^\infty}}{(t^*-s-3)^2}\,,
\end{aligned}
$$
for each $s\in (0,\tfrac12)$. Letting $s\to 0^+$ in the last inequality, one arrives
at
\begin{equation}
\label{limsup}
\begin{aligned}
\limsup_{\eps\to 0^+}\frac1{|\ln\eps|}\int_{\eps}^{\eps^*}\int_0^{\pi/4}
    \phi_r\left(\frac{t^*}r,(\cos\th,\sin\th)\right)m(\th)d\th\frac{dr}{r}&
\\
\le
\Lambda(t^*)+\frac{8\|m\|_{L^\infty}}{(t^*-3)^2}&\,.
\end{aligned}
\end{equation}
A similar inequality holds for the $\liminf$. By symmetry, the averaging in $\th$ in
(\ref{limsup}) can be done equivalently on $(0,\tfrac{\pi}4)$ or in $(0,2\pi)$, which
eventually proves (\ref{lim-dist}).

\section{Applications to kinetic theory}

It has been proved in Theorem 2.1 of \cite{GW} (see also \cite{BGW}) that the
linear Boltzmann equation (ie. equation (10) of \cite{Lo}) does {\it not} govern
the Boltzmann-Grad limit of the {\it periodic} Lorentz gas --- unlike the case
of a Lorentz gas with a {\it random} (Poisson) distribution of scatterers, where
the linear Boltzmann equation was rigorously derived in \cite{Ga} and \cite{BBS}.
While Theorem 2.1 of \cite{GW} is merely a negative result, we show below how
to infer from Proposition \ref{LIM-DIST} positive information on the asymptotic
behavior of the periodic Lorentz gas in the Boltzmann-Grad limit.

Define $\O_\eps=\{\eps z\,|\,z\in Z_\eps\}$, and consider the transport equation
\begin{equation}
\label{TranspEq}
\begin{aligned}
\d_tf_\eps+v\cdot\nabla_xf_\eps&=0\,,\quad x\in\O_\eps\,,\,\,|v|=1\,,
\\
f_\eps(t,x,v)&=0\,,\quad x\in\d\O_\eps\,,\,\,v\cdot n_x>0\,,
\\
f_\eps(0,x,v)&=f^{in}(x,v)\quad x\in\O_\eps\,,\,\,|v|=1\,.
\end{aligned}
\end{equation}
Here, the unknown is $f_\eps\equiv f_\eps(t,x,v)$ while $n_x$ is the inward unit
normal at point $x\in\d\O_\eps$ and $f^{in}$ is a given, nonnegative function of
$C_c(\R^2\times\S^1)$. Physically, this is a variant of the periodic Lorentz gas
where scatterers are replaced by holes (or traps) where impinging particles fall
and thus are removed from the domain $\O_\eps$. Obviously, for each $t\ge 0$
\begin{equation}
\label{MaxPr}
\|f_\eps\|_{L^\infty_{t,x,v}}=\|f^{in}\|_{L^\infty_{x,v}}\,.
\end{equation}
Reasoning as in \cite{Lo} suggests that $f_\eps\to f$ in $L^\infty_{t,x,v}$ weak-*
where $f$ solves the uniformly damped transport equation
\begin{equation}
\label{WrgKinLim}
\d_tf+v\cdot\nabla_xf+f=0\hbox{ on }\R_+^*\times\R^2\times\S^1\,,
    \quad f_{|t=0}=f^{in}\,,
\end{equation}
but this is ruled out by Theorem 2.1 of \cite{GW}. Instead, Proposition \ref{LIM-DIST}
suggests that the resulting damping rate should vanish in the limit as $t\to+\infty$.
This statement is made precise in the following theorem.

\begin{Thm}\label{LIM-KIN}
Let $f^{in}\ge 0$ belong to $C_c(\R^2\times\S^1)$ and let $f_\eps$ be, for each
$\eps\in(0,\tfrac14)$, the solution of (\ref{TranspEq}). Then, for each nonnegative
test function $\chi\in C^1_c(\R^2\times\S^1)$, one has
$$
\begin{aligned}
\limsup_{\eps\to 0}\iint\left(\frac1{|\ln\eps|}\int_{\eps}^{1/4}
    f_r(t,x,v)\frac{dr}{r}\right)\chi(x,v)dxdv&
\\
    =\iint f(t,x,v)\chi(x,v)dxdv+O\left(\frac1{t^2}\right)&
\\
\liminf_{\eps\to 0}\iint\left(\frac1{|\ln\eps|}\int_{\eps}^{1/4}
    f_r(t,x,v)\frac{dr}{r}\right)\chi(x,v)dxdv&
\\
    =\iint f(t,x,v)\chi(x,v)dxdv+O\left(\frac1{t^2}\right)&
\end{aligned}
$$
as $t\to+\infty$, where
\begin{equation}
\label{SolLim}
f(t,x,v)=\frac{2f^{in}(x-tv,v)}{\pi^2t}\,.
\end{equation}
In particular, $f$ satisfies
\begin{equation}
\label{KinLim}
\d_tf+v\cdot\nabla_xf+\frac{1}{t}f=0\,,
    \quad (t,x,v)\in(0,+\infty)\times\R^2\times\S^1
\end{equation}
in the sense of distributions.
\end{Thm}

\begin{proof}
First the solution of (\ref{TranspEq}) is given by the formula
\begin{equation}
\label{SolTransp}
f_\eps(t,x,v)=f^{in}(x-tv,v)\indc_{\tau_\eps(x/\eps,-v)\ge t/\eps}\,.
\end{equation}
In this formula, the exit time $\tau_\eps$ is considered as a function defined on
$Z_\eps\times\S^1$ with $\Z^2$-periodicity and extended by $0$ in $Z_\eps^c\times
\S^1$.

Let $h\equiv h(t,x,v)\ge 0$ belong to $C^\infty(\R_+\times\R^2\times\S^1)$, with
support in $\R_+\times[-L,L]^2\times\S^1$;
\begin{equation}
\label{kin-estim1}
\begin{aligned}
\left|\int\!\!\!\int h(t,x,v)\indc_{\tau_r(\frac{x}{r},-v)\ge\frac{t}r}dxdv
-
\int_{\S^1}\phi_r\left(\frac{t}r,v\right)\int h(t,x,v)dxdv\right|&
\\
\le
\left|\int\!\!\!\int h(t,x,v)\indc_{\tau_\eps(\frac{x}{r},-v)\ge\frac{t}r}dxdv
-
r^2\sum_{l\in\Z^2}\int_{\S^1}\!
    h(t,rl,v)\!\int_{\T^2}\!\indc_{\tau_r(y,-v)\ge\frac{t}r}dydv\right|&
\\
+
\int_{\S^1}\phi_r\left(\frac{t}r,v\right)\!\left|r^2\sum_{l\in\Z^2}h(t,rl,v)
    -\int h(t,x,v)dx\right|dv&
\\
=
r^2\sum_{l\in\Z^2}
    \iint_{\T^2\times\S^1}|h(t,rl+ry,v)-h(t,rl,v)|
    \indc_{\tau_r(y,-v)\ge\frac{t}r}dydv&
\\
+
r^2\sum_{l\in\Z^2}
    \iint_{\T^2\times\S^1}|h(t,rl+ry,v)-h(t,rl,v)|
    \phi_r\left(\frac{t}r,v\right)dydv&
\\
\le
2r\|\nabla_x h\|_{L^\infty}|S^1|\cdot r^2\sum_{l\in\Z^2}\indc_{[-L,L]^2}(rl)&
\\
\le
2r\|\nabla_x h\|_{L^\infty}|S^1|(L^2+o(1))&\,.
\end{aligned}
\end{equation}
By Theorem \ref{LIM-DIST}
\begin{equation}
\label{kin-estim2}
\begin{aligned}
\left|\limsup_{\eps\to 0^+}\frac1{|\ln\eps|}\int_\eps^{\eps^*}\int_{\S^1}
    \phi_r\left(\frac{t}r,v\right)\right.&\int h(t,x,v)dxdv\frac{dr}r
\\
&\left.-\frac2{\pi^2t}\iint h(t,x,v)dxdv\right|
\\
&\le
\frac{8\|h\|_{L^\infty_{t,v}(L^1_x)}}{t^*-3}
\end{aligned}
\end{equation}
with a similar estimate for the $\liminf$.

Putting together (\ref{kin-estim1}), (\ref{kin-estim2}) with the formula
(\ref{SolTransp}) establishes (\ref{SolLim}).
\end{proof}

Theorem \ref{LIM-KIN} can also be viewed as a result in homogenization. This remark
leads to a comparison with the analogous situation --- homogenization of a diffusion
process with Dirichlet boundary conditions in $Z_r$ --- studied in \cite{CiMu} (in
truth, the result obtained in \cite{CiMu} is much more complete and satisfying than
Theorem \ref{LIM-KIN}). Not surprisingly, the mathematical tools used in \cite{CiMu}
are of a very different nature than the ones in the present work. This however is by
no means surprising and simply reflects the very different nature of the trajectories
of a diffusion process and of those of a free transport equation. It is very likely
that the nature of the result in Theorem \ref{LIM-KIN} --- and its proof --- would be
deeply affected by adding some collision (ie. jump in velocity) process to the free
transport between successive impingements on the obstacles.

\section{Final remarks and perspectives}

A first question left open in the present work is the existence of the limit
(\ref{Q1}). As the reader will have probably noticed in the proof of Theorem
\ref{LIM-DIST}, the diameter $r$ of the obstacle is treated as a time variable
under renormalization --- ie. transformation of the problem by the Gauss map.
It may be that new insight concerning (\ref{Q1}) can be gained by using more
specific properties of $T$ than ergodicity as in the present paper.

A second question is the existence of the limit in the sense of Cesaro, as
commented upon in Remark \ref{CES-LIM}. As can be seen from the proof of
Theorem \ref{LIM-DIST}, proving this essentially amounts to being able to
apply an ergodic theorem as in section \ref{ERGO} to functions of the form
$$
f\left(\eps q_{N(\a,\eps)},
    \frac{d_N(\a,\eps)}{\eps},\frac{d_{N(\a,\eps)-1}}{\eps}\right)\,.
$$
(In the present paper, (\ref{Estim-qndn}) essentially allows one to replace
$q_{N(\a,\eps)}$ by $1/d_{N(\a,\eps)-1}$). This extension of the ergodic
theorem in the present paper is postponed to a subsequent paper \cite{CG}
and will be applied to the problem of the Lyapunov exponent for the Lorentz
gas (see \cite{Ch} for a presentation of this subject).

Finally we would like comment on some very interesting, related work in
\cite{Bo1} and \cite{Bo2}.

The billiard problem considered in \cite{Bo1} is similar to studying the
distribution of free path lengths in $Z_r$ with the obstacle at the origin
removed (ie. in $Z_r\cup\overline{B(0,r/2)}$) {\it for particles starting
from the origin only}. This is a quite different problem and thus the limit
as $r\to 0^+$ of this distribution has nothing to do with the simple result
in Theorem \ref{LIM-DIST}. A related issue is studied in \cite{Bo2}: for
$R>0$, consider the directions of lattice points (ie. point of $\Z^2$) in
the ball $B(0,R)$ of $\R^2$. These directions make a set of $N_R$ angles
$$
0\le\th_0<\th_1,\ldots<\th_{N_R}=2\pi\,;
$$
the paper \cite{Bo2} computes the distribution of scaled differences of the
form $N_R(\th_{k+1}-\th_k)$ as $R\to+\infty$. Remarkably, these differences
are far from being exponentially distributed --- unlike in the case of angles
picked at random in $[0,2\pi)$. Consistently with the results in \cite{BGW}
and in \cite{GW}, this observation stresses again the difference between the
case of a random distribution of scatterers as studied in \cite{Ga} and that
of a periodic distribution scatterers.

\section{Appendix:\\Background on continued fractions}\label{CONT-FR}

We recall below some basic facts and notations about continued fractions.
Given $\a\in(0,1)$, the Gauss map is defined by
\begin{equation}
\label{DefGauss}
T\a=\frac1{\a} -\left[\frac1{\a} \right]\,;
\end{equation}
it is known to be a uniquely ergodic transformation of $(0,1)$ with
invariant measure
\begin{equation}
\label{InvMeas}
dg(\a)=\frac1{\ln 2} \frac{d\a}{1+\a} \,.
\end{equation}
The continued fraction expansion of $\a$ is
\begin{equation}
\label{ContFrac}
\a=[a_1,a_2,a_3,\ldots]=
\frac1{a_1+{\displaystyle\frac1{a_2+{\displaystyle\frac1{a_3+\ldots}}}}}
    \hbox{ with }a_k=\left[\frac1{T^{k-1}\a} \right]\,,\,\,k\ge 1\,;
\end{equation}
and the action of $T$ is seen to correspond to the shift
\begin{equation}
\label{ActT}
\a=[a_1,a_2,a_3,\ldots]\mapsto T\a=[a_2,a_3,a_4,\ldots]\,.
\end{equation}
The convergents of $\a$ are defined by the recursion formulas
\begin{equation}
\label{DefCvgt}
\begin{aligned}
q_{n+1}&=a_nq_n+q_{n-1}\,,\quad q_0=0\,,\,\,q_1=1\,,
\\
p_{n+1}&=a_np_n+p_{n-1}\,,\quad p_0=1\,,\,\,p_1=0\,,
\end{aligned}
\end{equation}
and the corresponding error $d_n=|q_n\a-p_n|$ satisfies
\begin{equation}
\label{DefErr}
d_{n+1}=-a_nd_n+d_{n-1}\,,\quad d_0=1\,,\,\,d_1=\a\,.
\end{equation}
The first relation in (\ref{DefCvgt}) and (\ref{DefErr}) imply that
\begin{equation}
\label{qdRel}
q_nd_{n+1}+q_{n+1}d_n=1\,,\quad n\in\N^*\,.
\end{equation}
For each $\a\in(0,1)\setminus\Q$, the convergents of $\a$ are
the sequence of best rational approximants of $\a$. In other words
\begin{equation}
\label{BestAppr}
|q_n\a-p_n|
    =\inf\{\,|q\a-p|\,\,|\,\,p,q\in\Z\,,\,\,0\le q<q_n\}\,.
\end{equation}
This implies in particular that the approximation by continued
fractions is alternate in the sense that the algebraic value of
the error is
$$
q_n\a-p_n=(-1)^{n-1}d_n\,,\quad n\ge 0\,.
$$
Also, for each $\a\in(0,1)\setminus\Q$, the sequence of errors
satisfies the inequalities
\begin{equation}
\label{IneqErr}
\frac1{q_n+q_{n+1}}<d_n<\frac1{q_{n+1}}\,.
\end{equation}
The notations $a_n(\a)$, $p_n(\a)$, $q_n(\a)$ and $d_n(\a)$ are
used to emphasize the dependence of these quantities on $\a$
whenever necessary.

Our discussion uses ``renormalization'', ie. transforming $\a$ by
the iterates of the Gauss map $T$. We have gathered some useful
facts in the next lemma.

\begin{Lem}\label{LEM-REN}
Let $\a\in(0,1)$; then
\begin{itemize}
\item for each $n\in\N$, $q_n(T\a)=p_{n+1}(\a)$
    and $d_{n+1}(\a)=\a d_n(T\a)$;
\item for each $n\in\N$,
$$
d_n(\a)=\prod_{k=0}^{n-1}T^k\a\,;
$$
\end{itemize}
\end{Lem}

\begin{proof}
First observe that
$$
a_n(T\a)=a_{n+1}(\a)\,,\quad n\ge 1\,.
$$
This relation and (\ref{DefCvgt}) implies that the sequences $q_n(T\a)$
and $p_{n+1}(\a)$ satisfy the same recursion formulae; thus in order to
check the first formula, it suffices to check it for both $n=0$ and $n=1$
(in the latter case, one has $q_1(T\a)=[1/T\a]=a_2=p_2(\a)$). The second
formula is checked in the same way, and clearly implies the expression of
$d_n(\a)$ in terms of the $T^k\a$'s.
\end{proof}

\smallskip
\noindent
{\bf Acknowledgements.}
{\it We thank Prof. H. S. Dumas who told us that using ref. \cite{BK} might
simplify our original method for deriving formulas (\ref{Psi1'})-(\ref{Psi4'}).
The research of E. C. has been partially supported by MIUR and by INDAM GNFM.
Both authors acknowledge the support of the European Research Training Network
HyKE (Contract no. HPRN-CT-2002-00282).}

\end{document}